\documentclass{amsart}

\usepackage{amsmath,amssymb,amscd,amsfonts}
\usepackage[all]{xy}

\newtheorem{theorem}{Theorem}
\newcommand{\bt}{\begin{theorem}}
\newcommand{\et}{\end{theorem}}
\newtheorem{lemma}{Lemma}
\newcommand{\bl}{\begin{lemma}}
\newcommand{\el}{\end{lemma}}
\newtheorem{corollary}{Corollary}
\newcommand{\bc}{\begin{corollary}}
\newcommand{\ec}{\end{corollary}}
\newcommand{\beq}{\begin{equation}}
\newcommand{\eeq}{\end{equation}}
\newcommand{\benum}{\begin{enumerate}}
\newcommand{\eenum}{\end{enumerate}}
\newcommand{\N}{\ensuremath{ \mathbf N }}
\newcommand{\Z}{\ensuremath{\mathbf Z}}
\newcommand{\Q}{\ensuremath{\mathbf Q}}
\newcommand{\R}{\ensuremath{\mathbf R}}
\newcommand{\C}{\ensuremath{\mathbf C}}

\newcommand{\A}{\ensuremath{\mathbf A}}
\newcommand{\B}{\ensuremath{\mathbf B}}
\newcommand{\X}{\ensuremath{\mathbf X}}

\newcommand{\mbF}{\ensuremath{ \mathbf F}}

\newcommand{\mba}{\ensuremath{ \mathbf a}}
\newcommand{\mbb}{\ensuremath{ \mathbf b}}
\newcommand{\mbc}{\ensuremath{ \mathbf c}}

\newcommand{\mbx}{\ensuremath{ \mathbf x}}

\newcommand{\mbo}{\ensuremath{ \mathbf 0}}

\newcommand{\bmat}{\left(\begin{matrix}}
\newcommand{\emat}{\end{matrix}\right)}

\DeclareMathOperator{\vectorsmallan}{\left( \begin{smallmatrix} a_1 \\ \vdots \\ a_n \end{smallmatrix}\right)}
\DeclareMathOperator{\vectorsmallbn}{\left( \begin{smallmatrix} b_1 \\ \vdots \\ b_n \end{smallmatrix}\right)}
\DeclareMathOperator{\vectorsmallzero}{\left( \begin{smallmatrix} 0 \\ \vdots \\ 0 \end{smallmatrix} \right)}

\title{Sidon sets and perturbations}

\author{Melvyn B. Nathanson}
\address{Lehman College (CUNY),Bronx, New York 10468}
\email{melvyn.nathanson@lehman.cuny.edu}

\date{\today}



\thanks{Supported in part by a grant from the PSC-CUNY Research Award Program.}

\begin{document}

\maketitle

\begin{abstract}
A subset $A$ of an additive abelian group is an $h$-Sidon set if every element in the $h$-fold sumset 
$hA$ has a unique representation as the sum of $h$ not necessarily distinct elements of $A$.   
Let $\mathbf{F}$ be a field of characteristic 0 with a nontrivial absolute value, 
and let $A = \{a_i :i \in \mathbf{N} \}$ and $B = \{b_i :i \in \mathbf{N} \}$ be subsets of $\mathbf{F}$.
Let $\varepsilon =  \{  \varepsilon_i:i \in \mathbf{N} \}$,  where $\varepsilon_i > 0$ for all $i \in \mathbf{N}$.
The set $B$ is an $\varepsilon$-perturbation of  $A$ 
if $|b_i-a_i| < \varepsilon_i$ for all $i \in \mathbf{N}$.
It is proved that, for every $\varepsilon =   \{  \varepsilon_i:i \in \mathbf{N} \}$ with $\varepsilon_i > 0$,   
every set $A = \{a_i :i \in \mathbf{N} \}$  has an $\varepsilon$-perturbation $B$ 
that is an $h$-Sidon set.  This result extends to sets of vectors 
in $\mbF^n$.
\end{abstract}

\section{Sidon sets}
Let $\N\ = \{1,2,3,\ldots \}$ be the set of positive integers 
and $\N_0 = \N \cup \{ 0\} = \{0, 1, 2, 3, \ldots\}$ the set of nonnegative integers.  
Denote the cardinality of the set $S$ by $|S|$.  

Let  $G$ be an additive abelian group, 
and let  $A = \{a_i:i \in I \}$  be a nonempty subset of $G$ 
with $a_i \neq a_j$ for $i \neq j$.  
For every positive integer $h$, define the \emph{$h$-fold sumset}\index{sumset}
\[
hA = \{a_{i_1} + \cdots + a_{i_h}:  a_{i_1}, \ldots, a_{i_h} \in A \}.  
\] 
For every group element $b \in G$, we have $b \in hA$ if and only if 
there is a set of nonnegative integers $ \{ u_i: i \in I\}$ such that 
\[
h = \sum_{i \in I} u_i 
\]
and
\[
b = \sum_{i \in I} u_i a_i.
\]
Define  $0A = \{0\}$. 

An \emph{$h$-Sidon set}\index{Sidon set} (also called a \emph{$B_h$-set}) 
is a subset $A = \{a_i:i \in I \}$ of $G$ such that every element in the sumset 
$hA$ has a unique representation as the sum of $h$ elements of $A$.
Equivalently, an $h$-Sidon set is a set $A$ that satisfies the condition: 
For all $b \in hA$,  if $ \{ u_i: i \in I\}$ and $ \{ v_i: i \in I\}$ 
are sets of nonnegative integers such that 
\beq                    \label{perturb:uv-1}
h = \sum_{i \in I} u_i = \sum_{i \in I} v_i 
\eeq
and
\beq                    \label{perturb:uv-2}
b = \sum_{i \in I} u_i a_i = \sum_{i \in I} v_i a_i 
\eeq
then 
\beq                    \label{perturb:uv-3}
u_i = v_i \quad \text{ for all $i \in I$.}
\eeq

 A 2-Sidon set is usually called a \emph{Sidon set}.
 
A subset of an $h$-Sidon set is also an $h$-Sidon set, 
and an $h$-Sidon set is an $r$-Sidon set for every positive integer $r <  h$.  
However, an $h$-Sidon set is not necessarily an $(h+1)$-Sidon set.  
For example, let $G = \Z$ be the additive group of integers.  
For every integer $h \geq 2$, 
the set $A = \{ h^i: i \in \N_0 \}$ is an $h$-Sidon set 
because $0 \notin A$ and every positive integer has a unique $h$-adic representation. 
However, for all nonnegative integers $a < b$, the identity 
\[
h^{a+1} + h^{b+1}  = \underbrace{h^a + \cdots + h^a}_{\text{$h$ summands}}  + h^{b+1} = h^{a+1} + \underbrace{h^b + \cdots + h^b}_{\text{$h$ summands}} 
\]
proves that $A$ is not an $(h+1)$-Sidon set.

Most work on Sidon sets has been restricted to subsets of the integers 
or other discrete groups. 
Cilleruelo and Ruzsa~\cite{cill-ruzs04} have studied Sidon sets 
of real, complex, and $p$-adic numbers.  
O'Bryant~\cite{obry04} is a survey of Sidon sets.

\section{Perturbations of  countably infinite sets}

Let  $G$ be an additive abelian group and let $A$ and $B$ be subsets of $G$.  
Define the \index{difference set}\emph{difference set} 
\[
A - B =  \{a - b:  a\in A  \text{ and } b \in B \}   
\]
and, for all nonnegative integers $r$ and $s$,  the 
\index{sum-difference set}\emph{$(r,s)$-sum-difference set} 
\[
rA-sA = \left\{  \sum_{i=1}^r a_i - \sum_{i=r+1}^{r+s} a_i : 
a_i \in A \text{ for } i=1,\ldots, r+s  \right\}.
\]
For $c \in G$,  define the \emph{translate}\index{translate}  
\[
A+c = \{a+c: a \in A\}.
\] 
For $b \in G$ and $r \in \{0,1,\ldots, h\}$, let
\[
A_{r,h}(b) = rA + (h-r) b.  
\] 
We have 
\[
h \left(A  \cup \{ b \}\right) = \bigcup_{r =0}^h \left( rA + (h-r) b \right) 
= \bigcup_{r =0}^h A_{r,h}(b).
\]

\bl                                                                            \label{perturb:lemma:Arhb} 
Let  $G$ be an additive abelian group and let $A$ be a subset of $G$ 
that is an $h$-Sidon set.    Let $b \in G$.  The set $A \cup \{b\}$ is an $h$-Sidon set 
if and only if the  sets 
\[
A_{r,h}(b) = rA + (h-r) b  
\] 
are pairwise disjoint for all $r \in \{0,1,\ldots, h\}$.  
\el

\begin{proof}
Let $r,s \in \{0,1,\ldots, h\}$, and let $a_1,a_2,\ldots, a_{r+s}$ be a sequence of elements of $A$.  
We have  
\[
a_1 + \cdots + a_r + (h-r) b \in A_{r,h}(b)
\]
and 
\[
a_{r+1} + \cdots + a_{r+s} + (h-s) b \in A_{s,h}(b). 
\]
Suppose that 
\beq                                                                    \label{perturb:Arhb} 
c = a_1 + \cdots + a_r + (h-r) b = a_{r+1} + \cdots + a_{r+s} + (h-s) b.  
\eeq
Note that 
\[
c \in A_{r,h}(b) \cap A_{s,h}(b). 
\]
Because $A$ is an $h$-Sidon set, if $r=s$, then $a_1 + \cdots + a_r = a_{r+1} + \cdots + a_{2r}$ 
and there is a permutation 
$\sigma$ of $\{1,\ldots, r\}$ such that 
$a_{r+i} = a_{\sigma(i)}$ for all $i \in \{1,\ldots, r\}$.  

It follows that $A \cup \{b\}$ is not an $h$-Sidon set if and only if there exist integers 
$r,s \in \{0,1,\ldots, h\}$ with $r \neq s$, and elements $a_1,a_2,\ldots, a_{r+s}$ in $A$ 
that satisfy~\eqref{perturb:Arhb}.  This is equivalent to the condition 
that  the sets $A_{r,h}(b)$ are not pairwise disjoint for $r \in \{0,1,\ldots, h\}$.  
This completes the proof.  
\end{proof}

Let \mbF\ be a field.   For $A \subseteq \mbF$ and $c \in \mbF$,  define 
the \emph{dilate}\index{dilate}  
\[
c\ast A = \{ca:a \in A\}.
\]
An \index{absolute value} \emph{absolute value} $| \ |$ on the field \mbF\   
is trivial if $|0|=0$ and $|x| =1$ for all $x \in \mbF \setminus \{0\}$. 
A  field \mbF\ with a nontrivial absolute value is infinite, and 
\beq                       \label{perturb:AbsoluteValue} 
\inf\{ |x|:x \in \mbF\setminus \{0\} \} = 0.
\eeq 
The usual absolute values on \Q, \R, and \C\, and the $p$-adic absolute 
values on \Q\ and $\Q_p$ are nontrivial.  

Let \mbF\ be a field with a nontrivial absolute value $| \ |$. 
Let $I$ be a nonempty set, and let $\{a_i : i\in I \}$ and $\{b_i : i\in I \}$ 
be sets of elements of the field \mbF.  
Let $\{\varepsilon_i : i\in I \}$ be a set of positive real numbers.  
The set  $\{b_i : i\in I \}$  is an  
\emph{$\varepsilon$-perturbation}\index{perturbation}  
of the set $\{a_i : i\in I \}$ if 
\[
|b_i - a_i| < \varepsilon_i
\]
for all $i \in I$.

\bl                 \label{perturb:lemma:Sidon} 
Let \mbF\ be a field of characteristic 0 with a nontrivial absolute value $| \ |$, 
and let $A$ be a finite subset of \mbF.
For every $a \in \mbF$ and $\delta > 0$,   
there exists an element $b \in \mbF$ such that  
\[
| b-a | < \delta 
\]
and the $h+1$ sets 
\[
A_{r,h}(b) = rA + (h-r) b
\]
for all $ r \in \{0,1,2,\ldots, h\}$ are pairwise disjoint.
\el

\begin{proof}
 Let $a \in \mbF$ and $x \in \mbF$.   
 For all $r,s \in \{0,1,\ldots, h\}$ with $s < r$, we have 
\begin{align*} 
\emptyset & \neq 
A_{r,h}(a+x)   \bigcap A_{s,h}(a+x) \\
& = \left( rA + (h-r) (a+x)  \right) \bigcap \left(sA + (h-s) (a+x)  \right) 
\end{align*}
if and only if there is a sequence $a_1,a_2,\ldots, a_{r+s}$ of elements of $A$ such that 
\[
a_1 + \cdots + a_r + (h-r) (a+x)  = a_{r+1} + \cdots + a_{r+s} + (h-s) (a+x).   
\]
Equivalently,  
\beq                    \label{perturb:A}
(r-s) (a+x) = (a_1 + \cdots + a_r) - (a_{r+1} + \cdots + a_{r+s} ) \in rA - sA. 
\eeq
Because \mbF\ has characteristic 0, the positive integer $r-s$ is a unit in \mbF\ and 
\beq                    \label{perturb:B}
a+x \in  (r-s)^{-1} \ast (rA-sA).  
\eeq
Therefore, 
\[
x \in \left(  (r-s)^{-1}\ast (rA-sA) \right) - a.
\]
The set $A$ is finite, and so the sets 
\[
C_{r,s} =  \left(  (r-s)^{-1}\ast (rA-sA) \right) - a
\]
and 
\[
C = \bigcup_{0 \leq s < r \leq h} C_{r,s} 
\]
are also finite.   The set $\mbF \setminus C$ is infinite.  

For all $x \in \mbF \setminus C$, the sets $A_{r,h}(a + x)$ 
for $ r \in \{0,1,\ldots, h\}$ are pairwise disjoint.
Because the set $C$ is finite, we have  
\[
\delta_1 = \min\left(  |c|: c \in C \text{ and } c \neq 0 \right) > 0.
\]
The absolute value on \mbF\ is nontrivial, and so, by~\eqref{perturb:AbsoluteValue}, 
there exists $x\in \mbF$ 
with  $0 < |x| < \min(\delta_1,\delta)$.   
The inequality $0 < |x| < \delta_1$ implies that $x \notin C$, 
and so  the $h+1$ sets $A_{r,h}(a + x)$ are pairwise disjoint.
Let $b = a + x$.  We have $|b-a| =  |x|  < \delta$ and the sets 
$A_{r,h}(b) $ are pairwise disjoint for $ r \in \{0,1,\ldots, h\}$.
This completes the proof.  
\end{proof}

\bt                \label{perturb:theorem:FiniteSidon} 
Let \mbF\ be a field of characteristic 0 with a nontrivial absolute value $| \ |$. 
Let $A$ be a finite subset of \mbF\ that is an $h$-Sidon set.  
For every $a \in \mbF$ and $\delta > 0$,   
there exists an element $b \in \mbF$ such that  
\[
|b-a| < \delta 
\]
and $A \cup \{ b \}$  is an $h$-Sidon set.  
\et

\begin{proof}
This follows immediately from Lemmas~\ref{perturb:lemma:Arhb}  
and~\ref{perturb:lemma:Sidon}.  
\end{proof}

\bt                \label{perturb:theorem:InfiniteSidon} 
Let \mbF\ be a field of characteristic 0 with a nontrivial absolute value.   
Let $\{\varepsilon_i : i \in \N\}$ be a set of positive real numbers.  
For every subset $\{a_i: i \in \N \}$ of \mbF, there is an $h$-Sidon 
set $\{b_i: i \in \N \}$ in \mbF\  that is an $\epsilon$-perturbation of $A$.
\et

\begin{proof}
We construct the set $B$ by induction.  Let $b_1 = a_1$.
Let $k \geq 1$, and let $\{b_1,\ldots, b_k\}$ be an $h$-Sidon set in \mbF\ 
such that $|b_i - a_i| < \varepsilon_i$ for all $i \in \{1,\ldots, k\}$.  
Applying Theorem~\ref{perturb:theorem:FiniteSidon} with $A = \{b_1,\ldots, b_k\}$ 
and $a = a_{k+1}$, we obtain an element $b_{k+1} \in \mbF$ 
such that $| b_{k+1}  - a_{k+1} | < \varepsilon_{k+1}$ and 
$\{b_1,\ldots, b_k, b_{k+1}\}$ an $h$-Sidon set.
This completes the proof.  
\end{proof}

\bt           \label{perturb:corollary:InfiniteLimitSidon} 
Let \mbF\ be a field with a nontrivial absolute value. 
For every set  $A = \{a_i : i \in \N \}$ in \mbF\ 
there exists an $h$-Sidon set $B = \{ b_i : i \in \N \}$ in \mbF\ 
such that $\lim_{i\rightarrow \infty} |a_i - b_i| = 0$. 
\et

\begin{proof}
Apply Theorem~\ref{perturb:theorem:InfiniteSidon}  with a set 
$\{\varepsilon_i : i \in \N\}$ such that 
$\lim_{i\rightarrow \infty} \varepsilon_i = 0$.   
\end{proof}

\section{$h$-Sidon sets of vectors in $\mbF^n$} 

Let \mbF\ be a field with a nontrivial absolute value $| \ |$.   
For vectors 
$\mba = \vectorsmallan \in \mbF^n$ and $\mbb = \vectorsmallbn \in \mbF^n$, 
define 
\[
\|\mba - \mbb\| = \sum_{j =1}^n |a_j - b_j|. 
\]
Let $\mbo = \vectorsmallzero\in \mbF^n$ denote the zero vector.
It follows from~\eqref{perturb:AbsoluteValue} that 
\beq                       \label{perturb:VectorAbsoluteValue} 
\inf\{ \| \mbx\|: \mbx \in \mbF^n\setminus \{ \mbo\} \} = 0.
\eeq 
Let $c\mbx$ denote scalar multiplication of the vector \mbx\ by the scalar $c$.
Let $\X$ be a set of vectors.  The \emph{dilation of $\X$ by $c$} is the set  
$c\ast \X = \{c\mbx: \mbx \in \X\}$.
Note that if $h$ is a positive integer, then the dilated set $h \ast \X$ is a 
(usually proper) subset of the sumset $h\X$.  

\bl                 \label{perturb:lemma:VectorSidon} 
Let \mbF\ be a field of characteristic 0 with a nontrivial absolute value $| \ |$, 
and let $\A$ be a finite set of vectors in $\mbF^n$.
For every vector $\mba \in \mbF^n$ and $\delta > 0$,   
there exists a vector $\mbb \in \mbF^n$ such that  
\[
\| \mbb - \mba \| < \delta 
\]
and the $h+1$ sets 
\[
\A_{r,h}(\mbb) = r\A + (h-r) \mbb
\]
for all $ r \in \{0,1,2,\ldots, h\}$ are pairwise disjoint sets of vectors in $\mbF^n$.
\el

\begin{proof}
This is similar to the proof of Lemma~\ref{perturb:lemma:Sidon}.  
 Let $\mba \in \mbF^n$ and $\mbx \in \mbF^n$.   
 For all $r,s \in \{0,1,\ldots, h\}$ with $s < r$, we have 
\begin{align*} 
\emptyset & \neq 
\A_{r,h}(\mba + \mbx )   \bigcap \A_{s,h}(\mba + \mbx ) \\
& = \left( r\A + (h-r) (\mba + \mbx )  \right) \bigcap \left(s\A + (h-s) (\mba + \mbx )  \right) 
\end{align*}
if and only if there is a sequence $\mba_1,\mba_2,\ldots, \mba_{r+s}$ of vectors in $\A$ such that 
\[
\mba_1 + \cdots + \mba_r + (h-r) (\mba + \mbx )  = \mba_{r+1} + \cdots + \mba_{r+s} + (h-s) (\mba + \mbx ).   
\]
Equivalently,  
\[
(r-s) (\mba + \mbx ) = (\mba_1 + \cdots + \mba_r) - (\mba_{r+1} + \cdots + \mba_{r+s} ) \in r\A - s\A. 
\]
Because \mbF\ has characteristic 0, the positive integer $r-s$ is a unit in \mbF\ and 
\[
\mba + \mbx  \in  (r-s)^{-1} \ast (r\A-s\A).  
\]
Therefore, 
\[
\mbx \in \left(  (r-s)^{-1}\ast (r\A-s\A) \right) - \mba.
\]
The set $\A$ is finite, and so the sets 
\[
\C_{r,s} =  \left(  (r-s)^{-1}\ast (r\A-s\A) \right) - \mba
\]
and 
\[
\C = \bigcup_{0 \leq s < r \leq h} \C_{r,s} 
\]
are also finite.   The set $\mbF^n \setminus \C$ is infinite.  

For all vectors $\mbx \in \mbF^n \setminus \C$, the sets $\A_{r,h}(\mba + \mbx)$ 
for $ r \in \{0,1,\ldots, h\}$ are pairwise disjoint.
Because the set $\C$ is finite, we have  
\[
\delta_1 = \min\left(  |\mbc|: \mbc \in \C \text{ and } \mbc \neq \mbo \right) > 0.
\]
The absolute value on $\mbF^n$ is nontrivial, and so, by~\eqref{perturb:VectorAbsoluteValue}, 
there exists $\mbx\in \mbF^n$ 
with  $0 < \|\mbx \| < \min(\delta_1,\delta)$.   
The inequality $0 < \|\mbx \| < \delta_1$ implies that $\mbx \notin \C$, 
and so  the $h+1$ sets $\A_{r,h}(\mba + \mbx)$ are pairwise disjoint.
Let $\mbb = \mba + \mbx$.  We have $|\mbb - \mba| =  |\mbx|  < \delta$ and the sets 
$\A_{r,h}(\mbb) $ are pairwise disjoint for $ r \in \{0,1,\ldots, h\}$.
This completes the proof.  
\end{proof}

\bt                \label{perturb:theorem:FiniteSidon} 
Let \mbF\ be a field of characteristic 0 with a nontrivial absolute value $| \ |$. 
Let $\A$ be a finite set of vectors in $\mbF^n$ that is an $h$-Sidon set.  
For every vector $\mba \in \mbF^n$ and $\delta > 0$,   
there is a vector $\mbb \in \mbF^n$ such that  
\[
|\mbb - \mba| < \delta 
\]
and $\A \cup \{ \mbb \}$  is an $h$-Sidon set.  
\et

\begin{proof}
This follows immediately from Lemmas~\ref{perturb:lemma:Arhb}  
and~\ref{perturb:lemma:VectorSidon}.  
\end{proof}

\bt                \label{perturb:theorem:InfiniteVectorSidon} 
Let \mbF\ be a field of characteristic 0 with a nontrivial absolute value.   
Let $\varepsilon = \{\varepsilon_i : i \in \N\}$ be a set of positive real numbers.  
For every set $\A = \{\mba_i: i \in \N \}$ of vectors in $\mbF^n$, there is an $h$-Sidon 
set $\B = \{\mbb_i: i \in \N \}$ in $\mbF^n$  that is an $\varepsilon$-perturbation of $\A$.
\et

\begin{proof}
We construct the set $\B$ by induction.  Let $\mbb_1 = \mba_1$.
Let $k \geq 1$, and let $\{\mbb_1,\ldots, \mbb_k\}$ be an $h$-Sidon set in \mbF\ 
such that $|\mbb_i - \mba_i| < \varepsilon_i$ for all $i \in \{1,\ldots, k\}$.  
Applying Theorem~\ref{perturb:theorem:FiniteSidon} with $\A = \{\mbb_1,\ldots, \mbb_k\}$ 
and $\mba = \mba_{k+1}$, we obtain an element $\mbb_{k+1} \in \mbF$ 
such that $| \mbb_{k+1}  - \mba_{k+1} | < \varepsilon_{k+1}$ and 
$\{\mbb_1,\ldots, \mbb_k, \mbb_{k+1}\}$ an $h$-Sidon set.
This completes the proof.  
\end{proof}

\bt           \label{perturb:corollary:InfiniteLimitVectorSidon} 
Let \mbF\ be a field with a nontrivial absolute value. 
For every set of vectors $\A = \{\mba_i : i \in \N \}$ in $\mbF^n$
there is an $h$-Sidon set $\B = \{ \mbb_i : i \in \N \}$ in $\mbF^n$ 
such that $\lim_{i\rightarrow \infty} |\mba_i - \mbb_i| = 0$. 
\et

\begin{proof}
Apply Theorem~\ref{perturb:theorem:InfiniteSidon}  with a set 
$\{\varepsilon_i : i \in \N\}$ such that 
$\lim_{i\rightarrow \infty} \varepsilon_i = 0$.   
\end{proof}

\section{Open problems}
\benum
\item
The condition that the field \mbF\ has characteristic 0 is used only 
to deduce~\eqref{perturb:B}  from~\eqref{perturb:A}.
Let $A$ be a countably infinite set in a field of characteristic $p > 0$ with 
a nontrivial absolute value.  
Is there an $\varepsilon$-perturbation of $A$ that is an $h$-Sidon set? 
\item
Let $A$ be an uncountably infinite set.  Is there an $\varepsilon$-perturbation of $A$ that is an $h$-Sidon set? 
\eenum

\def\cprime{$'$} \def\cprime{$'$} \def\cprime{$'$}
\providecommand{\bysame}{\leavevmode\hbox to3em{\hrulefill}\thinspace}
\providecommand{\MR}{\relax\ifhmode\unskip\space\fi MR }
\providecommand{\MRhref}[2]{%
  \href{http://www.ams.org/mathscinet-getitem?mr=#1}{#2}
}
\providecommand{\href}[2]{#2}

\end{document}